\documentclass{svmult}

\usepackage{latexsym}

\newcommand{\RR}{\mathbf{R}}
\newcommand{\NN}{\mathbf{N}}

\newcommand{\F}{\mbox{$\mathcal{F}$}} 
\newcommand{\G}{\mbox{$\mathcal{G}$}}

\newcommand{\tN}{\mbox{$\tilde{N}$}} 
\newcommand{\tM}{\mbox{$\tilde{M}$}}

\newtheorem{Df}{Definition}[section]
\newtheorem{Th}[Df]{Theorem}

\newtheorem{Com}[Df]{Comment}

\begin{document}
\title*{A Compensator Characterization of Planar Point Processes}
\author{B. Gail Ivanoff\thanks{Research supported by a grant from the
Natural Sciences and Engineering Research Council of
Canada.}\\Dept. of Mathematics and Statistics\\University of Ottawa
\\585 King Edward Ave.\\Ottawa ON K1N 6N5 CANADA\\givanoff@uottawa.ca}
%\date{\today}
\authorrunning{Gail Ivanoff}
\maketitle
\vspace{-.5in}
\begin{center} {\em This paper is dedicated to Professor Mikl\'{o}s Cs\"{o}rg\H{o}, a wonderful mentor and friend, on the occasion of his $80^{th}$ birthday.}
\end{center}
\begin{abstract}\  
 Martingale techniques play a fundamental role in the analysis of point processes on $[0,\infty)$. In particular, the compensator of a point process uniquely determines and is determined by its distribution, and an explicit formula involving conditional interarrival distributions is well-known.   In two dimensions there are many possible  definitions of a point process compensator and we  focus here on the one that has been the most useful in practice: the so-called *-compensator.   Although existence of the *-compensator is well understood, in general it does not determine the law of the point process and it must be calculated on a case-by-case basis.
However, it will be proven that when the point process satisfies a certain property of conditional independence (usually denoted by (F4)), the *- compensator determines the law of the point process and an explicit regenerative formula can be given. The basic building block of the planar model is the {\em single line} process (a point process with incomparable jump points). Its law can be characterized by a class of avoidance probabilities that are the two-dimensional counterpart of the survival function on
$[0,\infty)$.  Conditional avoidance probabilities then play the same role in the construction of the *-compensator as conditional survival probabilities do for compensators in one dimension. 

%These facts have many applications, including a new martingale approach to nonparametric  inference for renewal processes in both one and two dimensions.
\end{abstract}
{\bf Keywords and phrases:}\\
point process, compensator, martingale, single line process, avoidance probability, cumulative hazard, adapted random set 
 \\ \\
{\bf AMS 2010 Subject Classifications} :  60G55, 60G48.
\\
\\
%{\bf Short Title}: Compensators for point processes.

\section{Background and Motivation}\label{section1}

If $N$ is a point process on $\RR_+$ with $E[N(t)]<\infty~\forall t\in \RR_+$, the compensator of $N$  is the unique predictable increasing process $\tN$ such that $N-\tN$ is a martingale  with respect to the minimal filtration generated by $N$, possibly augmented by information at time 0.   Why is $\tN$ so important? Some reasons include:
 \begin{itemize}
  \item  The law of $N$ determines and {\em is determined by} $\tN$ \cite{Jacod}. 
          \item  The asymptotic behaviour of a sequence of point processes can be determined by the asymptotic behaviour of the corresponding sequence of compensators \cite{{Brown}, {BIW}, {IMbook}}.
  \item   Martingale methods provide elegant and powerful nonparametric methods for point process inference, state estimation, change point problems, and easily incorporate censored data \cite{Karr}.     \end{itemize}

 Can martingale methods be applied to point processes in higher dimensions? This is an old question, dating back more than 30 years to the 70's and 80's when multiparameter martingale theory was an active area of research. However, since there are many different definitions of planar martingales, there is no single definition of ``the compensator" of a point process on $\RR_+^2$. A discussion of the various definitions and a more extensive literature review can be found in \cite{IMP} and \cite{IMbook}.
 
 In this article, we revisit  the following question:
 When can a compensator be defined for a planar point process in such a way that it exists, it is unique and it  characterizes the distribution of the point process? % Some related works include:
% Existence and uniqueness have been well studied in numerous contexts, including \cite{{Yor}, {MS}, {MZ}, {Dozzi}, {Gushchin}, {MM}, {Last}},%some relevant papers
 % \item Kallenberg  (1978): exvisible projections (not a martingale approach).
% \item Yor (1976), Mazziotto and Szpirglas (1980), 
%\item Last (1994):   existence of compensators for general point process filtrations. 
% but not the question of characterization. Some partial answers have been found:
%\begin{itemize}
%\item  In \cite{Jacod}, Jacod uses a marked point process approach to characterize the law of a planar point process by identifying one of the co-ordinates of each jump point as the mark of the jump.
%\item In \cite{MN}, compensator characterizations of the planar Poisson process are found. This is generalized in \cite{IMbook} to more general spaces.
%\item Merzbach and Zakai (1980), Dozzi (1981), Gushchin (1982): existence of compensators under (F4)
%\item  A compensator characterization of the law of a point process on $\RR_+^2$ is given for a specific two-dimensional filtration in \cite{IM90}.
%\item Ivanoff, Merzbach and Plante (2007):  general compensator characterization theorem via a family of one-dimensional compensators.% (not particularly useful!)
%\end{itemize}
Since there are many possible  definitions of a point process compensator in two dimensions, we  focus here on the one that has been the most useful in practice: the so-called *-compensator.   Although existence and uniqueness of the *-compensator is well understood \cite{Dozzi, Gushchin, Last}, in general it does not determine the law of the point process and it must be calculated on a case-by-case basis.
However, it will be proven in Theorem \ref{7.1} that when the point process satisfies a certain property of conditional independence (usually denoted by (F4), see  Definition \ref{F4}), the *- compensator determines the law of the point process and an explicit regenerative formula can be given. Although it seems to be widely conjectured that under (F4) the law must be characterized by the *-compensator, we have been unable to find a proof in the literature and, in particular, the related regenerative formula (\ref{8}) appears to be completely new.

 The basic building block of the planar model is the {\em single line} process (a point process with incomparable jump points). This approach was first introduced in \cite{MM} and further exploited in \cite{IMP}. In both cases, the planar process is embedded into a   point process with totally ordered jumps on a larger partially ordered space. ``Compensators" are then defined on the larger space. In the case of \cite{IMP},  this is   a family of one-dimensional  compensators that,  collectively, do in fact characterize the original distribution. Although the results in \cite{IMP} do not require the assumption (F4) and  are of theoretical significance, they seem to be difficult to apply in practice due to the abstract nature of the embedding.   So, although in some sense the problem of a compensator characterization has been  resolved for general planar point processes,  for practical purposes it is important to be able to work on the original space, $\RR_+^2$, if possible. We will see here that the assumption (F4) allows to do so.

Returning to the single line process,  when (F4) is satisfied we will see that its law can be characterized by a class of avoidance probabilities that form the two-dimensional counterpart of the survival function of a single jump point on $[0,\infty)$.  Conditional avoidance probabilities then play the same role in the construction of the *-compensator as conditional survival distributions do for compensators in one dimension. For clarity and ease of exposition, we will be assuming throughout continuity of the so-called avoidance probabilities; this will automatically ensure the necessary predictability conditions and connects the avoidance probabilities and the *-compensator via a simple logarithmic formula. %As well, we will assume that the underlying filtration is minimal. 
The more technical issues  of discontinuous avoidance probabilities and other related problems will be dealt  with in a separate publication. We comment further on these points in the Conclusion.

 Our arguments involve careful manipulation of conditional expectations with respect to different $\sigma$-fields, making repeated use of the conditional independence assumption (F4). For a good review of conditional independence and its implications, we refer the reader to \cite{Kallenberg}.

We proceed as follows: in \S \ref{section2}, we begin with a brief review of the  point process compensator on $\RR_+$, including its heuristic interpretation and its regenerative formula. In \S \ref{section3} we define compensators for planar point processes. We discuss the geometry and decomposition of planar point processes into ``single line processes" in \S \ref{section4}, and  in \S\ref{section5} we show how the single line processes can be interpreted via stopping sets, the two-dimensional analogue of a stopping time. The compensator of the single line process is developed in \S \ref{section6} and combined with the decomposition of  \S \ref{section4}, this leads in \S \ref{section7} to the main result, Theorem \ref{7.1}, which gives an explicit regenerative formula for the compensator of a planar point process  that characterizes its distribution. We conclude with some directions for further research in \S \ref{section8}.

\section{A quick review of the compensator on $\RR_+$}\label{section2}

There are several equivalent characterizations of a point process on $\RR_+$, and we refer the reader to \cite{DVJ} or \cite{Karr} for details. For our purposes, given a complete probability space $(\Omega, \F, P)$, we interpret a simple point process $N$ to be a pure jump  stochastic process  on $\RR_+$ defined by 
\begin{equation}\label{1} N(t):=\sum_{i=1}^\infty I(\tau_i\leq t),\end{equation}
where $0<\tau_1<\tau_2<... $  is a strictly increasing sequence of random variables (the jump points of $N$).
Assume that $E[N(t)]<\infty$ for every $t\in\RR_+$.  Let  $\F(t)\equiv\F_0\vee\F^N(t)$, where $\F^N(t):=\sigma\{N(s):s\leq t\}$, suitably completed, and $\F_0$ can be interpreted as information available at time 0.  This is a right-continuous filtration on $\RR_+$ and without loss of generality we assume $\F=\F(\infty)$. The law of $N$ is determined by its finite dimensional distributions.

Since $N$ is non-decreasing, it is an integrable submartingale and so has a Doob-Meyer decomposition $N-\tN$ where $\tN$ is the unique $\F$-predictable increasing process such that $N-\tN$ is a martingale. Heuristically,
$$\tN (dt)\approx P(N(dt)=1\mid \F(t)).$$
More formally, for each $t$, 
\begin{equation}\label{2}\tN (t)=\lim_{ n\rightarrow \infty }\sum_{k=0}^{2^n-1} E\left[N\left(\frac{(k+1)t}{2^n}\right)-N\left(\frac{kt}{2^n}\right)\mid \F\left(\frac{kt}{2^n}\right)\right],\end{equation}
where convergence is in the weak $L^1$-topology.

We have the following examples:
 \begin{enumerate} 

\item If $N$ is a Poisson process with mean measure $\Gamma$ and if $\F\equiv \F^N$, then by independence of the increments of $N$, it is an immediate consequence of (\ref{2}) that   $\tN=\Gamma$.
\item Let $N$ be a Cox process (doubly stochastic Poisson process): given a realization $\Gamma$ of a random measure $\gamma$ on $\RR_+$, $N$ is (conditionally) a Poisson process with mean measure $\Gamma$. If $\F_0=\sigma \{\gamma\}$, then $\tN=\gamma$. We refer to $\gamma$ as the driving measure of the Cox process.
\item The single jump process: Suppose that $N$ has a single jump point $\tau_1$, a r.v. with continuous distribution $F$ and let $\F\equiv \F^N$. In this case \cite{DVJ, Karr}
\begin{equation}\label{3}\tN(t)=\int_0^t I(u\leq \tau_1) \frac{dF(u)}{1-F(u )}=\Lambda(t\wedge \tau_1),\end{equation}
where $\Lambda(t):=-\ln(1-F(t ))$ is the cumulative (or integrated) hazard of $F$. $F$ is determined by its hazard $\frac{dF(\cdot)}{1-F(\cdot )}$. The relationship $\Lambda(t)=-\ln P(N(t)=0)$ in equation (\ref{3})  will be seen to have a direct analogue in two dimensions.
 
\item The general simple point process: We note that the jump points $(\tau_i)$ are $
\F$-stopping times and so we define $\F(\tau_i):=\{F\in \F:F\cap \{\tau_i\leq t\}\in \F(t)~\forall t\}$. Assume that for every $n$, there exists a continuous regular version $F_n(\cdot|\F(\tau_{n-1}))$ of the  conditional distribution of $\tau_n$ given $\F(\tau_{n-1})$     (we define $\tau_0=0$). Then if $\Lambda_n\equiv -\ln(1-F_n)$, we have the following regenerative formula for the compensator (cf. \cite{DVJ}, Theorem 14.1.IV):
\begin{equation}\label{4}
\tN (t)=\sum_{n=1}^\infty \Lambda_n(t\wedge \tau_n)I(\tau_{n-1}<t).
\end{equation}
Let $Q=P|_{{\cal F}_0}$ (the restriction of $P$ to $\F_0$). Since there is a 1-1 correspondence between $F_n$ and $\Lambda_n$, together, $Q$ and $\tN$ characterize  the law of $N$ (\cite{Jacod}, Theorem 3.4). When $\F\equiv \F^N$ (i.e. $\F_0$ is trivial), $\tN$ characterizes the law of $N$.
\end{enumerate}
\begin{Com}{\rm Note that $\Lambda_n$ can be regarded as a random measure with support on $(\tau_{n-1},\infty)$. Of course, in general we do not need to assume that $F_n$ is continuous in order to define the compensator (cf. \cite{DVJ}). However, the logarithmic relation above between $\Lambda _n$ and $F_n$ holds only in the continuous case, and we will be making analogous continuity assumptions for planar point processes.}
\end{Com}

\section{Compensators on $\RR_+^2$}\label{section3}

 We begin with some notation:
 For $s=(s_1,s_2),t=(t_1,t_2)\in \RR_+^2$, 
 \begin{itemize}
\item $s\leq t\Leftrightarrow s_1\leq t_1$ and $s_2\leq t_2$
\item $s\ll t\Leftrightarrow s_1< t_1$ and $s_2< t_2$.
\end{itemize}
We let $A_t:=\{s\in\RR_+^2:s\leq t\}$ and $D_t:=\{s\in \RR_+^2:s_1\leq t_1 \mbox{ or } s_2\leq t_2\}$. A set $L\subseteq \RR_+^2$ is a lower layer if  for every $t\in \RR_+^2$, $t\in L\Leftrightarrow A_t\subseteq L$. In analogy to (\ref{1}), given a complete probability space $(\Omega, \F, P)$ and distinct $\RR_+^2$-valued random variables $\tau_1,\tau_2,...$ (the jump points), the point process $N$ is defined by
\begin{equation}\label{5}
N(t):=\sum_{i=1}^\infty I(\tau_i\leq t)=\sum_{i=1}^\infty I(\tau_i\in A_t).
\end{equation}
As pointed out in \cite{Karr}, in $\RR_+^2$ there is no unique ordering of the indices of the jump points. Now letting $\tau_i=(\tau_{i,1},\tau_{i,2})$, we assume that $P(\tau_{i,1}=\tau_{j,1} \mbox{ for some } i\neq j)=P(\tau_{i,2}=\tau_{j,2} \mbox{ for some } i\neq j)=0$ and that $P(\tau_{i,1}=0)=P(\tau_{i,2}=0)=0~\forall i$. In this case, we say that
  $N$ is a strictly simple point process on $\RR_+^2$ (i.e. there is at most one jump point on each vertical and horizontal line and there are no points on the axes). The law of $N$ is determined by its finite dimensional distributions: $$P(N(t_1)=k_1,...N(t_i)=k_i),i\geq 1, t_1,...,t_i\in \RR_+^2, k_1,...k_i\in {\bf Z}_+.$$
  
  For any lower layer $L$, define $$\F^N(L):=\sigma (N(t):t\in L)$$
and \begin{equation}\label{oops}\F(L)=\F_0\vee \F^N(L),\end{equation} where $\F_0$ denotes the sigma-field of events known at time (0,0). In particular, since there are no jumps on the axes, $\F(L)=\F_0$ for $L$ equal to the axes. Furthermore, for any two lower layers $L_1$, $L_2$ it is easy to see that 
$$ \F(L_1)\vee \F(L_2) =\F(L_1\cup L_2) \mbox{ and }  \F(L_1)\cap \F(L_2) =\F(L_1\cap L_2).$$

 For $t\in \RR_+^2$, denote
$$\F(t):=\F(A_t)\mbox{ and } \F^*(t):=\F(D_t).$$
 Both $(\F(t))$ and $(\F^*(t))$ are right continuous filtrations indexed by $\RR_+^2$: i.e.  $\F^{(*)}(s)\subseteq \F^{(*)}(t)$ for all $s\leq t\in \RR_+^2$ and if  $t_n\downarrow t$,  then $\F^{(*)}(t)=\cap_n\F^{(*)}(t_n)$. More generally, if $(L_n)$ is a decreasing sequence of closed lower layers, $\F(\cap_n L_n)=\cap_n \F (L_n)$ (cf. \cite{IMMarkov}).

 \begin{Df}\label{martingale} Let $(X(t):t\in \RR_+^2)$  be an integrable stochastic process on $\RR_+^2$ and let $(\F(t):t\in \RR_+^2)$ be any filtration to which $X$ is adapted (i.e. $X(t)$ is $\F(t)$-measurable for all $t\in \RR_+^2$). $X$ is a weak $\F$-martingale if for any $s\leq t$, 
 $$E[X(s,t]\mid \F(s)]=0$$
 where $X(s,t]:=X(t_1,t_2)-X(s_1,t_2)-X(t_1,s_2)+X(s_1,s_2)$.
 \end{Df}
% \begin{picture}(200,180)(0,0)
%%\put(0,0){\vector(1,0){200}}
%\put(0,0){\vector(0,1){200}}
%\put(50,100){\circle*{5}}
%%\put(280,150){\circle*{5}}
%%\put(50,350){\circle*{5}}
%\put(25,40){$A_t$}
 %\put(50,100){\line(-1,0){50}}
% \put(50,100){\line(0,-1){100}}
%%\put(50,100){\line( 1,0){150}}
%%\put(50,100){\line(0,1){100}}
%\put(-10,-10){(0,0)}
%\put(50,-10){$t_1$}
%\put(-10,100){$t_2$}
%\put(53,103){$t$}
%\end{picture}

%  \begin{picture}(200,180)(0,0)
%\put(0,0){\vector(1,0){200}}
%\put(0,0){\vector(0,1){200}}
%\put(50,100){\circle*{5}}
%%\put(280,150){\circle*{5}}
%%\put(50,350){\circle*{5}}
%\put(25,40){$D_t$}
%% \put(50,100){\line(-1,0){50}}
%% \put(50,100){\line(0,-1){100}}
%\put(50,100){\line( 1,0){150}}
%\put(50,100){\line(0,1){100}}
%\put(-10,-10){(0,0)}
%%\put(50,-10){$t_1$}
%%\put(-10,100){$t_2$}
%\put(53,103){$t$}
%\end{picture}
% \\
% \\
 
We now turn our attention to point  process compensators on $\RR_+^2$. It will always be assumed that $E[N(t)]<\infty$ for every $t\in \RR_+^2$.
 For $t=(t_1,t_2)\in \RR_+^2$ and $0\leq k,j\leq 2^n-1$ define
 \begin{eqnarray*} \Delta N(k,j)&:= &
 N\left(\left(\frac{kt_1}{2^n},\frac{jt_2}{2^n}\right),\left(\frac{(k+1)t_1}{2^n},\frac{(j+1)t_2}{2^n}\right)
 \right]%\\
% &=&N\left(\frac{(k+1)t_1}{2^n},\frac{(j+1)t_2}{2^n}\right)-N\left(\frac{(k+1)t_1}{2^n},\frac{jt_2}{2^n}\right)\\
 %&&-N\left(\frac{kt_1}{2^n},\frac{(j+1)t_2}{2^n}\right)+N\left(\frac{kt_1}{2^n},\frac{jt_2}{2^n}\right)
 .
 \end{eqnarray*}
 %oAssume that $N$ has finite mean.
 In analogy to $\RR_+$, the {\em weak $\F$}-compensator of $N$ is defined by
 $$\tN (t):=\lim_{ n\rightarrow \infty }\sum_{j=0}^{2^n-1}\sum_{k=0}^{2^n-1} E\left[\Delta N(k,j) \mid \F\left(\frac{kt_1}{2^n},\frac{jt_2}{2^n}\right)\right],$$
and the $\F^*$-compensator ({\em strong} $\F$-compensator) of $N$ is defined by  
$$\tN^* (t):=\lim_{ n\rightarrow \infty }\sum_{j=0}^{2^n-1}\sum_{k=0}^{2^n-1} E\left[\Delta N(k,j) \mid \F^*\left(\frac{kt_1}{2^n},\frac{jt_2}{2^n}\right)\right],$$
where both limits are in the weak $L^1$ topology. When there is no ambiguity, reference to $\F$ will be suppressed in the notation. Note that although $\tN^*$ is $\F^*$-adapted, it is not $\F$-adapted in general.

\begin{Com} {\rm
  Under very general conditions, the compensators  exist and     $N-\tN$ and $N-\tN^*$ are  weak  martingales with respect to $\F$ and $\F^*$, respectively \cite{IMbook, Last}.   Furthermore, each has a type of predictability property that ensures uniqueness (cf.  \cite{IMbook}). Both compensators have non-negative increments: $\tN^{(*)}(s,t]\geq 0~\forall s,t\in \RR_+^2$.
   However, neither compensator determines the distribution of $N$ in general, as can be seen in the following examples.  }
 \end{Com}
   
\noindent{\bf Examples:}
  \begin{enumerate}
 \item The Poisson and Cox processes: Let $N$ be a Poisson process on $\RR_+^2$ with mean measure $\Gamma$ and let $\F=\F^N$. By independence of the increments, both the weak and *-compensators  of $N$ ($\tN$ and $\tN^*$) are equal to $\Gamma$ (\cite{IMbook}, Theorem 4.5.2). A deterministic *-compensator characterizes the Poisson process, but a deterministic weak compensator does not (see \cite{IMbook} for details). Likewise, if $N$ is a Cox process with driving measure $\gamma$ on $\RR_+^2$ and if $\F_0=\sigma\{\gamma\}$, then $\tN^*\equiv\gamma$; this too characterizes the Cox process (cf. \cite{IMbook}, Theorem 5.3.1). This discussion can be summarized as follows:
 \begin{Th}\label{3.3} Let $N$ be a strictly simple point process on $\RR_+^2$ and let $\gamma$ be a random measure on $\RR_+^2$ that puts 0 mass on every vertical and horizontal line. Let  $\F_0=\sigma\{\gamma\}$ and $\F(t)=\F_0\vee \F^N(t),~ \forall t\in \RR_+^2$. Then $N$ is a Cox process with driving measure $\gamma$ if and only if $\tN^*\equiv \gamma$. The law of $N$ is therefore determined by $Q:=P|_{{\cal F}_0}$ and $\tN^*$. In the case that $\gamma$ is deterministic, $\F_0$ is trivial and $N$ is a Poisson process.
 \end{Th}
 \item  The  single jump process: Assume that $N$ has a single jump point $\tau\in \RR_+^2$, a random variable with continuous distribution $F$   and survival function $S(u)=P(\tau\geq u)$. Then (cf. \cite{IMbook}):%(\neq 1-F)$. 
\begin{eqnarray*}
\tN (t)&=&\int_{[0,{t_1}]\times[0,{t_2}]} I(u\leq \tau) \frac{dF(u)}{1-F(u)}, \mbox{ and}\\
\tN^*(t)&=&\int_{[0,{t_1}]\times[0,{t_2}]} I(u\leq \tau) \frac{dF(u)}{S(u)}. 
\end{eqnarray*}
Although both formulas look very similar to (\ref{3}), in two dimensions  it is well known that neither $dF(u)/(1-F(u))$ nor $dF(u)/S(u)$ determines $F$.  \end{enumerate}

% weak and * compensator - non uniqueness and comments -when does the compensator determine the %distribution\\
% * not adapted, Poisson, single jump\\
% What information is needed, and the BIG question
So we see that neither $\tN$ nor $\tN^*$ determines the law of $N$ in general. The problem is that the filtration $\F$ does not provide enough information about $N$, and in some sense the filtration $\F^*$ can provide too much. As was observed in \cite{IMP}, the correct amount of information at time $t$ lies between $\F(t)$ and $\F^*(t)$.  The solution would be to identify a condition under which  the two filtrations provide essentially the same information  -  this occurs under a type of conditional independence, a condition usually denoted by (F4) in the two-dimensional martingale literature. 

  To be precise, for  $t=(t_1,t_2)\in \RR_+^2$ and any filtration $(\F(t))$, define the following $\sigma$-fields:
  \begin {eqnarray*}%\F(t)&=&\sigma\{N(s):s\in A_t\}= \sigma\{N(s):s_1\leq t_1\mbox{ and } s_2\leq t_2\}
%\label{Fdef} \\
 \F^1(t)&:=&\vee_{s\in \RR_+}\F(t_1,s)\nonumber\\
 \F^2(t)&:=&\vee_{s\in \RR_+}\F(s,t_2) \nonumber. 
 \end{eqnarray*} 
\begin{Df}\label{F4}
We say that the filtration $(\F(t))$ satisfies condition $(F4)$ if for all $t\in\RR_+^2$, the $\sigma$-fields $\F^1(t)$ and $\F^2(t)$ are conditionally independent, given $\F(t)$ ($\F^1(t)\perp \F^2(t)\;|\; \F(t)$).
 \end{Df}
 
 For the point process filtration $\F(t)=\F_0\vee \F^N(t)$, in practical terms (F4) means that the behaviour of the point process is determined only by points in the past (in terms of the partial order): geographically, this means by points from the southwest.  $N$ could denote the  points of infection in the spread of an air-born disease under prevailing winds from the southwest: since there are no points in $[0,t_1]\times(t_2,\infty)$ southwest of $(t_1,\infty) \times[0,t_2]$ and vice versa, the behaviour of $N$ in either region will not affect the other.
 
 While it appears that (F4) is related to the choice of the axes, it can be expressed in terms of the partial order on ${\RR_+^2}$. In fact, it is equivalent to the requirement that for any $s, t\in \RR_+^2$, 
 $$E[E[\cdot\mid \F(s)]\mid\F(t)] =E[\cdot \mid \F(s\wedge t).]]$$
 This concept can be extended in a natural way to other  partially ordered spaces; see Definition 1.4.2 of  \cite{IMbook}, for example. 

% {\bf Comment:} Strictly speaking, it is the law of $N$ that satisfies (F4).  
  % \begin{picture}(200,180)(0,0)
%\put(0,0){\vector(1,0){200}}
%\put(0,0){\vector(0,1){200}}
%\put(50,100){\circle*{5}}
%\put(280,150){\circle*{5}}
%\put(50,350){\circle*{5}}
%\put(25,40){$\F(t)$}
%\put(18, 140){$\F^1(t)$}
%\put(125, 40){$\F^2(t)$}
% \put(50,100){\line(-1,0){50}}
% \put(50,100){\line(0,-1){100}}
%\put(50,100){\line( 1,0){150}}
%\put(50,100){\line(0,1){100}}
%\put(-10,-10){(0,0)}
%\put(50,-10){$t_1$}
%\put(-10,100){$t_2$}
%\put(53,103){$t$}
%\end{picture}

%Note: the single jump process does not satisfy (F4).

%SECOND TRY
 Condition (F4) has the following important consequence: if $F\in \F(t)=\F_0\vee\F^N(t)$, then for any lower layer $D$, 
  \begin{equation}\label{add}P[F\mid\F(D)]=P[F\mid\F(t)\cap \F(D)].  \end{equation}
  This is proven in \cite{Dozzi} for $D=D_s$ for  $s\in \RR_+^2$, and the result is easily generalized as follows. To avoid trivialities, assume $t\not\in D$. Let $s_1:=\sup\{s\in \RR_+:(s,t_2)\in D\}$ and $s_2:=\sup\{s\in \RR_+:(t_1,s)\in D\}$ and define the lower layers $D_1$ and $D_2$ as follows:
  \begin{eqnarray*}
  D_1&:=&\{u=(u_1,u_2)\in D:u_1\leq s_1\}\\
   D_2&:=&\{u=(u_1,u_2)\in D:u_2\leq s_2\}
 \end{eqnarray*}%
 We have that $D=(D\cap A_t)\cup D_1\cup D_2$ and $\F(D)=\F(A_t\cap D)\vee \F(D_1)\vee\F(D_2)$.   By (F4), $ \F(D_2)\perp (\F(t)\vee \F(D_1))\;|\: \F((t_1,s_2))$. Now use the chain rule for conditional expectation (\cite{Kallenberg}, Theorem 5.8):
  \begin{eqnarray}
 \lefteqn{ \F(D_2)\perp (\F(t)\vee \F(D_1))\;|\: \F((t_1,s_2))}\nonumber\\
 &\Rightarrow & \F(D_2)\perp \F(t)\;|\: (\F((t_1,s_2))\vee \F(D_1))\nonumber\\
  &\Rightarrow & \F(D_2)\perp \F(t)\;|\: (\F((t_1,s_2))\vee \F(D_1)\vee \F(A_t\cap D))\label{l1}\\
 &\Rightarrow & \F(D_2)\perp \F(t)\;|\: ( \F(D_1)\vee \F(A_t\cap D)).\label{l2} 
 \end{eqnarray}
(\ref{l1}) and (\ref{l2}) follow  since $\F((t_1,s_2))\subseteq \F(A_t\cap D)\subseteq \F(t)$. But once again by (F4) we have $\F(D_1)\perp \F(t)\;|\; \F((s_1,t_2))$, and since $\F((s_1,t_2))\subseteq \F(A_t\cap D)\subseteq \F(t)$ we have
\begin{eqnarray}
\F(D_1)\perp \F(t)\;|\; \F((s_1,t_2))&\Rightarrow& \F(D_1)\perp \F(t)\;|\;( \F((s_1,t_2))\vee \F(A_t\cap D))\nonumber\\
&\Rightarrow &\F(D_1)\perp \F(t)\;|\;  \F({A_t\cap D}). \label{l3} 
\end{eqnarray}
Finally, if $F\in \F(t)$, 
\begin{eqnarray*}
P[F\mid \F_D]&=&P[F\mid \F({A_t\cap D})\vee \F({D_1})\vee\F({D_2})]\\
&=& P[F\mid \F(A_t\cap D)\vee \F(D_1)] \mbox{ by (\ref{l2})}\\
&=&P[F\mid \F(A_t\cap D)] \mbox{ by (\ref{l3})},
\end{eqnarray*}
and (\ref{add}) follows since $ \F({A_t\cap D})=\F(t) \cap \F(D)$.

We can use (\ref{add}) to argue heuristically that  (F4) ensures   that $\F$ and $\F^*$ provide roughly the same information: 
 \begin{eqnarray*}\lefteqn{E\left[\Delta N(k,j) \mid \F^*\left(\frac{kt_1}{2^n},\frac{jt_2}{2^n}\right)\right]  }\\
& =& E\left[\Delta N(k,j) \mid \F \left(\frac{(k+1)t_1}{2^n},\frac{jt_2}{2^n}\right)\vee \F\left(\frac{kt_1}{2^n},\frac{(j+1)t_2}{2^n}\right)\right]  \mbox{by (F4) (cf. (\ref{add})) } \\
 & \approx & E\left[\Delta N(k,j) \mid \F \left(\frac{kt_1}{2^n},\frac{jt_2}{2^n}\right)\right]
  \mbox{ as $n\rightarrow \infty$}.
  \end{eqnarray*}
  Therefore, $ \tN\approx \tN^*$  and in particular, $\tN^*$ is $\F$-adapted. In this case, we refer to  $N-\tN^*$ as a {\em strong} $\F$-martingale:
   \begin{Df}\label{smartingale} Let $(X(t):t\in \RR_+^2)$  be an integrable stochastic process on $\RR_+^2$ and let $(\F(t):t\in \RR_+^2)$ be any filtration to which $X$ is adapted. $X$ is a strong $\F$-martingale if for any $s\leq t$, 
 $$E[X(s,t]\mid \F^*(s)]=0.$$
 \end{Df}  
  
   As mentioned before, to avoid a lengthy discussion of predictability we will deal only with continuous compensators. In this case, we have the following (cf. \cite{Dozzi,Gushchin}):
  
  \begin{Th}\label{strongmg} Let $N$ be a strictly simple point process and assume that the filtration $\F=\F_0\vee \F^N$ satisfies (F4). If $\gamma$ is a continuous increasing $\F$-adapted process such that $N-\gamma$ is a strong martingale, then $\tN^*\equiv \gamma$. (We say that $\gamma$ is increasing if $\gamma(s,t]\geq 0~\forall~s\leq t\in \RR_+^2$.)

  \end{Th}
  
We now address the following question: if (F4) is satisfied, will the *-compensator characterize the distribution of $N$? In the case of both the Poisson and Cox processes, (F4) is satisfied for the appropriate filtration ($\F(t)=\F^N(t)$ for the Poisson process and $\F(t)=\sigma\{\gamma\}\vee \F^N(t)$ for the Cox process) and the answer is yes, as noted in Theorem \ref{3.3}. For these two special cases, it is possible to exploit one dimensional techniques since  conditioned on $\F_0$, the *-compensator is deterministic (see \cite{IMbook}, Theorem 5.3.1). Unfortunately, this one dimensional approach cannot be used for more general point process compensators. Nonetheless, Theorem \ref{3.3} turns out to be the key to the general construction of the compensator. 

Before continuing, we note here that  when (F4) is assumed a priori and the point process is strictly simple, 
there are many other characterizations of the two-dimensional Poisson process - for a thorough discussion see \cite{Merzbach}.  Assuming (F4), another approach  is to project the two-dimensional point process onto a family of increasing paths. Under different sets of  conditions, it is shown in \cite{AC} and \cite{IMP} that if the compensators of the corresponding one-dimensional point processes are deterministic, the original point process is Poisson. (For a comparison of these results, see \cite{IMP}.) However, the characterization of the Poisson and Cox processes given in Theorem \ref{3.3} does {\em not} require the hypothesis of (F4), and in fact implies it. Furthermore, it can be extended to more general spaces and to  point processes that are not strictly simple (cf. \cite{IMbook}, Theorem 5.3.1), although (F4) will no longer necessarily be satisfied.

Returning to the general case, the first step is to analyze the geometry of strictly simple point processes from the point of view taken in   \cite{IMP} and \cite{MM}.
  
 \section{The geometry of point processes on $\RR_+^2$}\label{section4}
    
Let $d=1$ or 2. If  $N$ is a strictly simple point process  on $\RR_+^d$,  then $N$ can be characterized via the increasing family of random sets 
 $$\xi_k(N):=\{t\in \RR_+^d: N(s)<k ~\forall s\ll t\}, k\geq 1.$$
 By convention,   in $\RR_+$ we define $\xi_0(N)$ to be the origin, and in  $\RR_+^2 $ we define $\xi_0(N)$ to be  the axes.
  We observe that:
 \begin{itemize}
\item In $\RR_+$, $\xi_k(N)=[0,\tau_k]$. 
\item  $N(t)=k\Leftrightarrow t\in\{\xi_{k+1}^o(N)\setminus\xi_k^o(N)\}.$ ($\xi_k^o(N)$ denotes the interior of $\xi_k(N)$.)
\item  In $\RR_+^2$, $\xi_k(N)$ is defined by the set of its {\em exposed points}:
 $${\cal E}_k:=\min\{t\in \RR_+^2:N(t)\geq k\}$$
 where for a nonempty Borel set $B\subseteq \RR_+^d$, $\min(B):=\{t\in B: s\not\leq t,\forall s\in B, s\neq t\}$. By convention, $\min(\emptyset):= \infty$.
 It is easily seen that 
$$ \xi_k(N)=\cap_{\tau \in {\cal E}_k} D_\tau.$$ 
 \end{itemize}
 
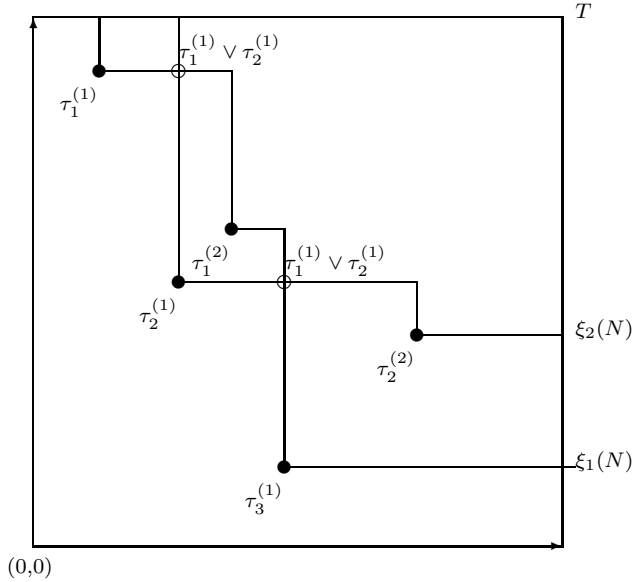
\begin{figure}  
\begin{center} 
\begin{picture}(200,200)(0,0) 
\put(0,0){\vector(1,0){200}} 
\put(0,0){\vector(0,1){200}} 
 \put(55,100){\circle*{5}} 
%\put(55,100){\dashbox{.8}(40,80)} 
 \put(25,180){\circle*{5}}
%\put(25,180){\dashbox{.8}(30,20)}  
 \put(95,30){\circle*{5}}
% \put(95,30){\dashbox{.8}(105,70)}
 %\put(175,150){\circle*{5}} 
\put(145,80){\circle*{5}} 
\put(75,120){\circle*{5}} 
%\put(75,120){\dashbox{.5}(20,60)} 
\put(95,30){{\line(1,0){105}}} 
 \put(95,30){{\line(0,1){70}}} 
\put(55,100){{\line(1,0){40}}} 
 \put(55,100){{\line(0,1){80}}}  
 \put(25,180){{\line(1,0){30}}} 
 \put(25,180){{\line(0,1){20}}} 
%\put(95,180){{\line(0,-1){60}}} 
%\put(55,180){\dashbox{.5}(20,20)} 
%\put(145,80){\dashbox{.5}(55,20)} 
%\put(25,180){\line(1,0){175}} 
%\put(25,180){\line(0,1){20}} 
%\put(95,30){\line(1,0){105}} 
%\put(95,30){\line(0,1){170}} 
%\put(175,150){\line(1,0){25}} 
%\put(175,150){\line(0,1){50}} 
\put(145,80){\line(1,0){55}}
\put(145,80){\line(0,1){20}} 
\put(95,100){\line(1,0){50}}
\put(95,30){\line(0,1){90}}%changed
\put(95,30){\line(1,0){110}}%changed
\put(75,120){\line(1,0){20}}
\put(75,120){\line(0,1){60}}
\put(55,180){\line(1,0){20}}
\put(55,180){\line(0,1){20}}
%\put(145,80){\line(0,1){120}} 
%\put(75,120){\line(1,0){125}} 
%\put(75,120){\line(0,1){80}} 
%\put(75,180){\circle{5}} 
%\put(95,180){\circle{5}} 
%\put(145,180){\circle{5}} 
%\put(175,180){\circle{5}} 
 \put(55,180){\circle{5}} 
 \put(95,100){\circle{5}} 
%\put(95,100){\dashbox{.5}(50,20)} 
%\put(95,120){\line(1,0){50}} 
%\put(145,100){\line(1,0){55}} 
%\put(145,120){\line(0,-1){20}} 
%\put(95,120){\circle{5}} 
%\put(145,100){\circle{5}} 
%\put(145,120){\circle{5}} 
\put(-10,-10){(0,0)} 
\put(0,0){\framebox(200,200)} 
\put(205,200){$T$} 
  \put(40,85){$\tau^{(1)}_2$} 
  \put(10,165){$\tau^{(1)}_1$} 
  \put(80,15){$\tau^{(1)}_3$} 
%\put(160,135){$\tau^{(5)}_2$} 
  \put(130,65){$\tau^{(2)}_2$} 
 \put(60,105){$\tau^{(2)}_1$} 
%\put(60,165){$\tau^{(3)}_1$} 
%\put(80,165){$\tau^{(4)}_1$} 
%\put(130,165){$\tau^{(5)}_1$} 
%\put(160,165){$\tau^{(6)}_1$} 
  \put(55,185){$\tau^{(1)}_1\vee\tau^{(1)}_2$} 
  \put(95,105){$\tau^{(1)}_1\vee\tau^{(1)}_2$} 
%\put(80,105){$\tau^{(3)}_2$} 
%\put(130,85){$\tau^{(3)}_3$} 
%\put(130,105){$\tau^{(4)}_2$} 
%\put(205,95){$\xi_1^+(N)$} 
\put(205,30){$\xi_1(N)$} 
\put(205,80){$\xi_2(N)$} 

% \put(25,190){\tiny$ E^o_{\tau^{(1)}_1}\cap \xi_1^+$} 
%\put(57,140){\tiny$ E^o_{\tau^{(1)}_2}\cap \xi_1^+$} 
%\put(130,60){\tiny$ E^o_{\tau^{(1)}_3}\cap \xi_1^+$} 
\end{picture} 
\end{center} 
\caption{Upper boundaries of the random sets $\xi_1(N)$  and $\xi_2(N)$. Jump points of $N$ indicated by {\large $\bullet$}. Other exposed points indicated by {\large$\circ $}. }
\label{Fig1}
\end{figure}

 To illustrate, in Figure \ref{Fig1} we consider the random sets $\xi_1(N)$  and $\xi_2(N)$ of a point process with five jump points, each indicated by a ``{\large $\bullet$}"'.   While the exposed points $\tau^{(1)}_1,\tau^{(1)}_2,\tau^{(1)}_3$ of $\xi_1(N)$ are all jump points of $N$, the exposed points of $\xi_2(N)$ include  $\tau^{(1)}_1\vee\tau^{(1)}_2$ and $\tau^{(1)}_2\vee\tau^{(1)}_3$ (each indicated by a ``{\large$\circ $}") , which are not jump points. In fact, if
 $$\xi_k^+(N):=\cap_{\epsilon,\epsilon'\in{\cal E}_k, \epsilon\neq \epsilon '}\ D_{\epsilon\vee\epsilon '},$$
then $\xi_k(N)\subseteq \xi_{k+1}(N)\subseteq \xi_k^+(N)$.
If ${\cal E}_k$ is empty or consists of a single point, then $\xi_k^+(N):=\RR_+^2$. For the same example, the upper boundaries of the sets  $\xi_1(N)$ and $\xi_1^+(N)$ are illustrated in Figure \ref{Fig2}.

%\begin{center} 
%\end{center}  
  
 \begin{figure}  
\begin{center} 
\begin{picture}(200,200)(0,0) 
\put(0,0){\vector(1,0){200}} 
\put(0,0){\vector(0,1){200}} 
 \put(55,100){\circle*{5}} 
\put(55,100){\dashbox{.8}(40,80)} 
 \put(25,180){\circle*{5}}
\put(25,180){\dashbox{.8}(30,20)}  
 \put(95,30){\circle*{5}}
 \put(95,30){\dashbox{.8}(105,70)}
 %\put(175,150){\circle*{5}} 
%\put(145,80){\circle*{5}} 
%\put(75,120){\circle*{5}} 
%\put(75,120){\dashbox{.5}(20,60)} 
\put(95,30){{\line(1,0){105}}} 
 \put(95,30){{\line(0,1){70}}} 
\put(55,100){{\line(1,0){40}}} 
 \put(55,100){{\line(0,1){80}}} 
 \put(25,180){{\line(1,0){30}}} 
 \put(25,180){{\line(0,1){20}}} 
%\put(95,180){{\line(0,-1){60}}} 
%\put(55,180){\dashbox{.5}(20,20)} 
%\put(145,80){\dashbox{.5}(55,20)} 
%\put(25,180){\line(1,0){175}} 
%\put(25,180){\line(0,1){20}} 
%\put(95,30){\line(1,0){105}} 
%\put(95,30){\line(0,1){170}} 
%\put(175,150){\line(1,0){25}} 
%\put(175,150){\line(0,1){50}} 
%\put(145,80){\line(1,0){55}}
%\put(145,80){\line(0,1){20}} 
%\put(95,100){\line(1,0){50}}
\put(95,30){\line(0,1){70}}%changed
\put(95,30){\line(1,0){110}}%changed
%\put(75,120){\line(1,0){20}}
%\put(75,120){\line(0,1){60}}
%\put(55,180){\line(1,0){20}}
%\put(55,180){\line(0,1){20}}
%\put(145,80){\line(0,1){120}} 
%\put(75,120){\line(1,0){125}} 
%\put(75,120){\line(0,1){80}} 
%\put(75,180){\circle{5}} 
%\put(95,180){\circle{5}} 
%\put(145,180){\circle{5}} 
%\put(175,180){\circle{5}} 
 \put(55,180){\circle{5}} 
 \put(95,100){\circle{5}} 
%\put(95,100){\dashbox{.5}(50,20)} 
%\put(95,120){\line(1,0){50}} 
%\put(145,100){\line(1,0){55}} 
%\put(145,120){\line(0,-1){20}} 
%\put(95,120){\circle{5}} 
%\put(145,100){\circle{5}} 
%\put(145,120){\circle{5}} 
\put(-10,-10){(0,0)} 
\put(0,0){\framebox(200,200)} 
\put(205,200){$T$} 
  \put(40,85){$\tau^{(1)}_2$} 
  \put(10,165){$\tau^{(1)}_1$} 
  \put(80,15){$\tau^{(1)}_3$} 
%\put(160,135){$\tau^{(5)}_2$} 
 % \put(130,65){$\tau^{(2)}_2$} 
 %\put(60,105){$\tau^{(2)}_1$} 
%\put(60,165){$\tau^{(3)}_1$} 
%\put(80,165){$\tau^{(4)}_1$} 
%\put(130,165){$\tau^{(5)}_1$} 
%\put(160,165){$\tau^{(6)}_1$} 
  \put(55,185){$\tau^{(1)}_1\vee\tau^{(1)}_2$} 
  \put(95,105){$\tau^{(1)}_2\vee\tau^{(1)}_3$} 
%\put(80,105){$\tau^{(3)}_2$} 
%\put(130,85){$\tau^{(3)}_3$} 
%\put(130,105){$\tau^{(4)}_2$} 
\put(205,95){$\xi_1^+(N)$} 
\put(205,30){$\xi_1(N)$} 
%\put(205,80){$\xi_2(N)$} 

% \put(25,190){\tiny$ E^o_{\tau^{(1)}_1}\cap \xi_1^+$} 
%\put(57,140){\tiny$ E^o_{\tau^{(1)}_2}\cap \xi_1^+$} 
%\put(130,60){\tiny$ E^o_{\tau^{(1)}_3}\cap \xi_1^+$} 
\end{picture} 
\end{center} 
 \caption{Upper boundaries of the random sets $\xi_1(N)$ and $ \xi_1^+(N)$.} 
\label{Fig2}
\end{figure}
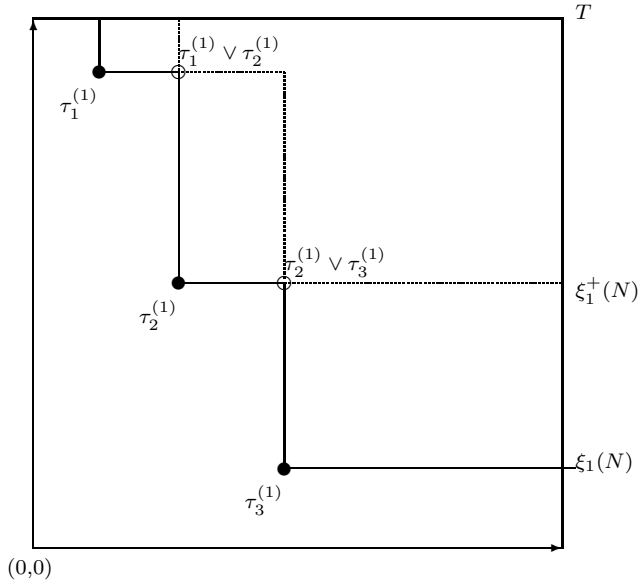
 
\begin{figure}  

\begin{center} 
\begin{picture}(200,200)(0,0) 
\put(0,0){\vector(1,0){200}} 
\put(0,0){\vector(0,1){200}} 
 \put(55,100){\circle*{5}} 
\put(55,100){\dashbox{.8}(40,80)} 
 \put(25,180){\circle*{5}}
\put(25,180){\dashbox{.8}(30,20)}  
 \put(95,30){\circle*{5}}
 \put(95,30){\dashbox{.8}(105,70)}
 %\put(175,150){\circle*{5}} 
\put(142,78){\large$\otimes$} 
\put(72,118){\large$\otimes$} 
%\put(75,120){\dashbox{.5}(20,60)} 
\put(95,30){{\line(1,0){105}}} 
 \put(95,30){{\line(0,1){70}}} 
\put(55,100){{\line(1,0){40}}} 
 \put(55,100){{\line(0,1){80}}} 
 \put(25,180){{\line(1,0){30}}} 
 \put(25,180){{\line(0,1){20}}} 
%\put(95,180){{\line(0,-1){60}}} 
%\put(55,180){\dashbox{.5}(20,20)} 
%\put(145,80){\dashbox{.5}(55,20)} 
%\put(25,180){\line(1,0){175}} 
%\put(25,180){\line(0,1){20}} 
%\put(95,30){\line(1,0){105}} 
%\put(95,30){\line(0,1){170}} 
%\put(175,150){\line(1,0){25}} 
%\put(175,150){\line(0,1){50}}
{ 
%\color{alert}
\put(145,80){\line(1,0){55}}
\put(145,80){\line(0,1){20}} 
\put(95,100){\line(1,0){50}}
\put(95,30){\line(0,1){90}}%changed
\put(95,30){\line(1,0){110}}%changed
\put(75,120){\line(1,0){20}}
\put(75,120){\line(0,1){60}}
\put(55,180){\line(1,0){20}}
\put(55,180){\line(0,1){20}}}
%\put(145,80){\line(0,1){120}} 
%\put(75,120){\line(1,0){125}} 
%\put(75,120){\line(0,1){80}} 
%\put(75,180){\circle{5}} 
%\put(95,180){\circle{5}} 
%\put(145,180){\circle{5}} 
%\put(175,180){\circle{5}} 
 \put(55,180){\circle{5}} 
 \put(95,100){\circle{5}} 
%\put(95,100){\dashbox{.5}(50,20)} 
%\put(95,120){\line(1,0){50}} 
%\put(145,100){\line(1,0){55}} 
%\put(145,120){\line(0,-1){20}} 
%\put(95,120){\circle{5}} 
%\put(145,100){\circle{5}} 
%\put(145,120){\circle{5}} 
\put(-10,-10){(0,0)} 
\put(0,0){\framebox(200,200)} 
\put(205,200){$T$} 
  \put(40,85){$\tau^{(1)}_2$} 
  \put(10,165){$\tau^{(1)}_1$} 
  \put(80,15){$\tau^{(1)}_3$} 
%\put(160,135){$\tau^{(5)}_2$} 
  \put(130,65){$\tau^{(2)}_2$} 
 \put(60,105){$\tau^{(2)}_1$} 
%\put(60,165){$\tau^{(3)}_1$} 
%\put(80,165){$\tau^{(4)}_1$} 
%\put(130,165){$\tau^{(5)}_1$} 
%\put(160,165){$\tau^{(6)}_1$} 
  \put(55,185){$\tau^{(1)}_1\vee\tau^{(1)}_2$} 
  \put(95,105){$\tau^{(1)}_2\vee\tau^{(1)}_3$} 
%\put(80,105){$\tau^{(3)}_2$} 
%\put(130,85){$\tau^{(3)}_3$} 
%\put(130,105){$\tau^{(4)}_2$} 
\put(205,95){$\xi_1^+(N)$} 
\put(205,30){$\xi_1(N)$} 
\put(205,80){$\xi_2(N)$} 

% \put(25,190){\tiny$ E^o_{\tau^{(1)}_1}\cap \xi_1^+$} 
%\put(57,140){\tiny$ E^o_{\tau^{(1)}_2}\cap \xi_1^+$} 
%\put(130,60){\tiny$ E^o_{\tau^{(1)}_3}\cap \xi_1^+$} 
\end{picture}

\end{center} 
\caption{Jump points of $M_1$ indicated by {\Large$\bullet$}, jump points of    $M_2$ indicated by{\large$\otimes$}.} 
\label{Fig3}
\end{figure}
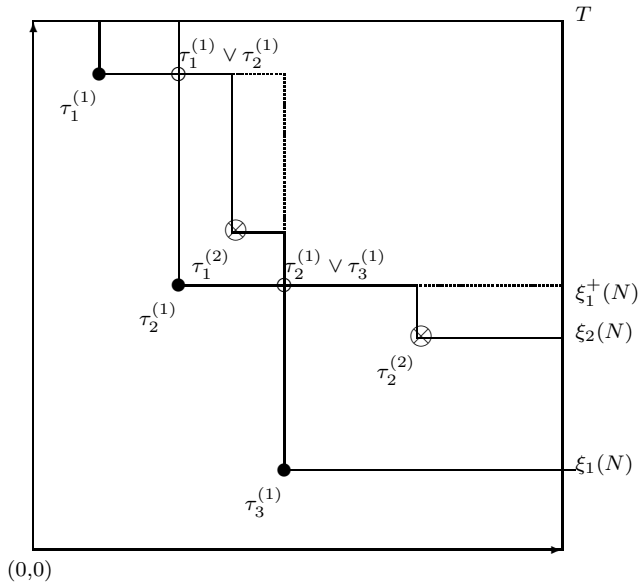
 
 We can now define $N$ in terms of {\em single line} point processes:  
\begin{Df} A point process on $\RR_+^2$ whose jump points are all   incomparable  is a single line process. (Points $s,t\in \RR_+^2$ are incomparable if $s\not\leq t$ and $t\not\leq s$.)\end{Df}
 \begin{Df}\label{4.2} Let $N$ be a strictly simple point process on $\RR_+^2$ and let $J(N)$ denote the set of jump points of $N$. Then $N(t)=\sum_1^\infty M_k(t)$ where for $k\geq 1$, $M_k$ is the single line process whose set of jump points is  $$J(M_k):=\min\left(J(N)\cap (\xi_{k-1}^+(N)\setminus\xi_{k-1}(N))\right),$$
where  $\xi_0=\{\{0\}\times \RR_+\}\cup \{\RR_+\times\{0\}\}$ and $\xi_0^+=\RR_+^2$.
  \end{Df}
   Returning to our example, in Figure \ref{Fig3} we illustrate each of the jump points of $M_1$ with {\Large$\bullet$}, and each of  jump points of    $M_2$  with {\large$\otimes$}. 
   
Before continuing, we make a few observations:
\begin{itemize} 
\item $\xi_{k}(N)=\xi_{k-1}^+(N)\cap \xi_1(M_k)$ (this is illustrated in Figure \ref{Fig3} for $k=2$). We note that $M_k$ has no jump points if $J(M_k)=\emptyset$; in this case $\xi_1(M_k)=\RR_+^2$ and $\xi_{k}(N)=\xi_{k-1}^+(N)$.
%\item The single line process is the natural two-dimensional analogue of the single jump process.
\item Since $\{N(t)=k\}=\{t\in \xi_{k+1}^o(N)\setminus \xi_{k}^o(N)\}$, in a manner that will be made precise, the law of $N$ (its finite dimensional distributions) is determined by the joint (finite dimensional) distributions of the random sets $\xi_k^o(N)$. %(equivalently, $\xi_1(M_k)$). 
We will see that this can be done by successive conditioning,  as in one dimension where the joint distribution of the successive jump times is built up through conditioning. 
\item If $M$ is a single line process, it is completely determined by $\xi_1(M)$ (cf. \cite{IMP} - the jump points of $M$ are the exposed points of $\xi_1(M)$). 
\item  Since the point process and its related random sets $\xi_k(N)$ are determined by single line processes, we will be able to reduce our problem  to the following question:  will the *-compensator of the single line process $M_k$ characterize its distribution if (F4) is satisfied? 
 \end{itemize}
 First, we need to consider the concept of stopping in higher dimensions.
 
 \section{Stopping sets and their distributions}\label{section5}
 
 We begin with the definition of adapted random sets and stopping sets; in particular, a stopping set is the multidimensional analogue of a stopping time. %First, $L\subseteq \RR_+^d$ is a lower layer if $t\in L\Rightarrow A_t\subseteq L~\forall t \in\RR_+^d$.  
 \begin{Df} Let $d=1$ or $2$. An adapted random set $\zeta$ with respect to the  filtration $\F$ on $\RR_+^d$  is a random Borel subset of $\RR_+^d$ such that $\{t\in \zeta\}\in \F(t)~\forall t\in \RR_+^d$. An adapted random set $\xi$ is an $\F$-stopping set if $\xi$ is a closed lower layer.
 \end{Df}
 For $d=1$, we see that if $\tau$ is an $\F$-stopping time, then $\zeta=[0,\tau)$ is an adapted random set and $\xi=[0,\tau]$ is an $\F$-stopping set. Since $\F$ is right-continuous, it is easily seen that $\xi=[0,\tau]$ is an $\F$-stopping set if and only if $\tau$ is an $\F$-stopping time. For $d=2$ and $\F(t)=\F_0\vee \F^N(t)$ for a point process $N$, and if $\F(L)$ is  defined as in (\ref{oops}) for a lower layer $L$, then it is shown in \cite{IMbook} that both $\{\xi\subseteq L\}\in \F(L)$ and $\{L\subseteq \xi\}\in \F(L)$.
 
The law of an adapted random set $\zeta$ is determined by its finite dimensional distributions:  $$P(t_1,...,t_n\in \zeta), n\in\NN, t_1,...,t_n\in\RR_+^d, d=1 \mbox{ or } 2.$$ 

%Given a filtration $\{\F(t):t\in \RR_+^2\}$ and $L$ a closed lower layer, define $\F(L):=\vee_{t\in L}\F(t)$. (For a point process $N$, this definition coincides with (\ref{oops}).)
In analogy to the history of a stopping time, the history of  a stopping set $\xi$ is 
 $$\F(\xi):=\{G\in \F: G\cap \{\xi\subseteq L\}\in \F(L)~ \forall\mbox{ lower layers } L\}.$$ If $\xi$ takes on at most countably many values in the class of lower layers, then equality can be used in the definition above, and it is easy to see that  $\F(\xi)=\F(L)$ on $\{\xi=L\}$.
For any point process $N$ on $\RR_+^2$ and filtration $\F\equiv\F_0\vee \F^N$,  we have the following:
\begin{itemize}
\item  Since $\{t\in \xi_k^o\}=\{N(t)<k\}\in \F_t~\forall k, t$, $\xi_k^o(N)$ is an $\F$-adapted random set. \item It is shown in \cite{IMrenewal} that the sets $\xi_k(N)$ and $\xi_k^+(N)$ are both $\F$-stopping sets. As well,   both  are $\F(\xi_k(N))$-measurable for every $k$ (i.e. $\{t\in\xi_k^{(+)}(N)\}   \in \F(\xi_k)~\forall t$).
\item Since $N$ is strictly simple, a priori there are no jumps on the axes and so $\F(\xi_0)=\F_0$. 
\end{itemize} 

Just as the joint distributions of the increasing jump times $\tau_1<\tau_2<...$ determine the law of a point process on $\RR_+$ and can be built up by successive conditioning on $\F_0\subseteq\F(\tau_1)\subseteq \F(\tau_2)\subseteq...$, we see that the law (finite dimensional distributions) of a planar point process $N$ can be reconstructed from the joint finite dimensional distributions of the related adapted random sets:
$$P(N(t_1)=k_1,...N(t_n)=k_n)=P(t_i\in \xi_{k_i+1}^o\setminus \xi_{k_i}^o,i=1,...,n).$$
As well, it is clear that the joint distributions of the increasing random sets $\xi_1^o\subset \xi_2^o\subset ..$ can be built up by successive conditioning on $\F_0=\F(\xi_0)\subseteq\F(\xi_1)\subseteq\F(\xi_2)\subseteq...$ .

 \section{The compensator of a single line point process}\label{section6} 
 
We are now ready to construct the  *-compensator of a single line process $M$ on $\RR_+^2$. Of course, we continue to assume that $E[M(t)]<\infty~\forall t\in \RR_+^2$. 

Although in principle the law of a point process is determined by the joint laws of the sets $\xi_k^o(M), k\geq 1$, in the case of a single line process, the law of $M$ is completely determined by the law of $\xi_1^o(M)$ (\cite{IMP}, Proposition 5.1). In other words, the set of probabilities
$$P(M(t_1)=0,...,M(t_n)=0)=P(t_1,...,t_n\in \xi_1^o(M))$$
for $t_1,...,t_n\in \RR_+^2$, $n\geq 1$, characterize the law of $M$. (This can be compared with the characterization of the law of a point process on an arbitrary complete measurable metric space via the so-called avoidance function; see \cite{DVJ}, Theorem 7.3.II.) 

However, when (F4) is satisfied, we have a further simplification. Define the {\em avoidance probability function} $P_0$ of a single line process $M$ by $$P_0(t):=P(M(t)=0), t\in \RR_+^2.$$
%\begin{Lem}\label{avoidance}(\cite{IMP}, Lemma 5.3)m
%\end{Lem} Let $M$ be a single line process on $\RR_+^2$ satisfying (F4). The law of $M$ is completely determined by the {\em avoidance probability function} $P_0$, defined by
%$$P_0(t)=P(M(t)=0)=P(t\in \xi_1^o(M), ~ t\in \RR_+^2).$$ 
 \begin{Th}\label{6.1}(\cite{IMP}, Lemma 5.3) Let $M$ be a single line process whose minimal filtration $\F\equiv\F^M$ satisfies (F4). The law (the f.d.d.'s) of $\xi_1^o(M)$ (and hence the law of $M$) is determined by the avoidance probability function $P_0$ of $M$. 
    \end{Th} 

A complete proof is given in \cite{IMP}, but to illustrate, we consider two incomparable points $s,t\in \RR_+^2$. If $s_1<t_1$ and $t_2< s_2$,  recalling that $\F=\F^M$ satisfies  (F4) and   that $M$ is a single line process, we have: 
\begin{eqnarray*}\lefteqn{
P(M(t)=0\mid\F^1(s))}\\&=&P(M(t)=0\mid\F (s\wedge t))  \mbox{ by (F4) (cf. (\ref{add}))}\\
&=&P(M(t)=0\mid   M(s\wedge t)=0 )I(M(s\wedge t)=0) \\
&=&\frac{P(M(t)=0) }{P(M(s\wedge t)=0 )}I(M(s\wedge t)=0).
\end{eqnarray*}
Therefore, 
\begin{eqnarray}\lefteqn{ 
P(M(s)=0,M(t)=0) }\nonumber\\&=& P(s,t\in \xi_1^o(M)) \nonumber \\
%&=& P(M(s\wedge t)=0,P(M((s\wedge t),t]=0), P(M((s\wedge t),s]=0)\\
&=& E\left[I(M(s)=0) P(M(t)=0\mid\F^1(s))\right]\nonumber \\
&=& E\left[I(M(s)= 0)I(M(s\wedge t)=0)\frac{P(M(t)=0)}{P(M(s\wedge t)=0)}\right]\nonumber \\
&=& E\left[I(M(s)= 0)\frac{P(M(t)=0)}{P(M(s\wedge t)=0)}\right]\nonumber\\
&=& \frac{P(M(s)=0)P(M(t)=0)}{P(M(s\wedge t)=0)}=\frac{P_0(s)P_0(t)}{P_0(s\wedge t)}.\label{6}
\end{eqnarray}

%In fact, letting $\{N(t-)=0\}:=\{N(s)=0~\forall s\llt\}$ 
Under (F4), the  avoidance probability function of a single line process can be regarded as the two-dimensional analogue of the survival function of the jump time $\tau$ of a single jump point process on $\RR_+$. {\em Henceforth, we will assume that the avoidance probability function is continuous.}
Obviously, the avoidance probability function is non-increasing in the partial order on $\RR_+^2$, but when is a continuous function bounded by 0 and 1 and non-increasing in each variable an avoidance probability? The answer lies in its logarithm. 

Let $\Lambda (t):=-\ln P_0(t)=-\ln P(M(t)=0)$.  Returning to (\ref{6})  and taking logarithms on both sides, if $s,t\in \RR_+^2$ are incomparable,
\begin{eqnarray*}
\Lambda(s\vee t)& =&-\ln P(M(s\vee t)=0)\\
&\geq& -\ln P(M(A_s\cup A_t)=0)\mbox{ since $A_s\cup A_t\subseteq A_{s\vee t}$}\\
&=&-\ln P(M(s)=0,M(t)=0)\\
&=& \Lambda ( s)+\Lambda( t)-\Lambda (s\wedge t) \mbox{ by (\ref{6})}.
\end{eqnarray*}
 If $P_0$ is continuous, then $\Lambda$ is continuous and increasing on $\RR_+^2$: i.e. it has non-negative increments. Therefore, $\Lambda=-\ln P_0$ is the distribution function of  a measure on $\RR_+^2$. In what follows, we will use the same notation for both the measure and its distribution function; for example, for $B$ a Borel set, $\Lambda (B)$ and $M(B)$ are the measures assigned to $B$ by the distribution functions $\Lambda(t)=\Lambda(A_t)$ and $M(t)=M(A_t)$, respectively. To summarize, when $\F=\F^M$ satisfies (F4):   % Return to the question: when will the compensator of a single line process  determine  its distribution?

  \begin{itemize}
  \item { If $P_0$ is continuous, $\Lambda =-\ln P_0$ defines a measure on $\RR_+^2$, and it is straightforward that for any lower layer $L$,}
  $$P(L\subseteq   \xi_1(M))=e^{-\Lambda(L)}=P(M(L)=0).$$ 
 \item{ Conversely, a  measure $\Lambda$ that puts mass 0 on each vertical and  horizontal line uniquely defines the (continuous) avoidance probability function $P_0$ (and therefore the law) of a single line point process whose minimal filtration  satisfies (F4).}
\item {Heuristically, $d\Lambda$ can be interpreted as the hazard of $M$:} 
$$P(M(dt)=1|\F^*(t))\stackrel{(F4)}{\approx}I(M(A_t)=0)d\Lambda(t).
$$ We will refer to $\Lambda$ as the {\em cumulative hazard} of $M$.
%\item All of the preceding discussion can be applied to conditional avoidance functions and conditional cumulative hazard functions. START HERE _ NEED TO POINT OUT THAT (F4) WORKS WITH F0 AS WELL.k
  \end{itemize}

All of the preceding discussion can be applied to conditional avoidance probability functions and conditional cumulative hazard functions, but first we need to define regularity of conditional avoidance probabilities; this is analogous to the definition of a regular conditional distribution.  
\begin{Df}\label{regular} Given an arbitrary $\sigma$-field  ${\cal F'}\subseteq \F$, we say that a family $(P_0(t ,\omega):(t,
\omega)\in \RR_+^2\times \Omega)$ is a continuous regular version of a conditional avoidance probability function given ${\cal F'}$ if for each $t\in \RR_+^2$, $P_0(t ,\cdot)$ is ${\cal F'}$-measurable, and for each $\omega\in \Omega$, $P_0(\cdot , \omega)$ is  equal to one on the axes, and $-\ln P_0(\cdot , \omega)$ is continuous and increasing on $\RR_+^2$.
\end{Df}
 We have the following generalization of Theorem \ref{6.1}:

 \begin{Th}\label{6.3} Let $M$ be a single line process with filtration $\F(t)= \F_0\vee\F^M(t)$ that satisfies (F4). If there exists a continuous regular version  $P_0^{(0)}(\cdot,\omega)$  of the conditional avoidance probability of $M$ given $\F_0$, then the conditional law   of $\xi_1^o(M)$ (and hence $M$) given $\F_0$ is determined by $P_0^{(0)}$, or equivalently by the conditional cumulative hazard $\Lambda_0:=-\ln  P_0^{(0)}$.    \end{Th}

Now we can define the *-compensator of the single line process; to do so, we will make use of Theorem \ref{3.3}. Suppose first that we have the minimal filtration: $\F(t)=
\F^M(t)$. Since $P(M(L)=0)=e^{-\Lambda(L)}$ for any lower layer $L$,
we can identify  the single line process $M$ with the single line process $M_1$ (the first line)  in the decomposition of a Poisson process $N$ with continuous mean measure $\Lambda$ (cf. Definition \ref{4.2}): we have $\xi_1(M)=\xi_1(M_1)=\xi_1(N)$. As shown in Example 7.4 of \cite{IMrenewal}, it is easy to see that the  $(\F^N)^*$-compensator of $M$ is $\tM^*(t)=\Lambda(A_t\cap \xi_1(M))$. However, since $\F^M\subseteq \F^N$ and $\tM^*$ is $\F^M$-adapted, by Theorem \ref{strongmg} it follows that $\tM^*$ is also the $(\F^M)^*$-compensator of $M$. Similarly, if 
$\F(t)=\F_0\vee
\F^M(t)$, since $P(M(L)=0\mid\F_0)=e^{-\Lambda_0(L)}$ for any lower layer $L$,  we make the same identification with a Cox process with driving measure $ \Lambda_0 $ to obtain  $\tM^*(t)=\Lambda_0(A_t\cap \xi_1(M))$ (as above, this is both the $(\F_0\vee\F^N)^*$ and the $(\F_0\vee\F^M)^*$-compensator).
We summarize this as follows:

    \begin{Th}\label{6.4} Let $M$ be a single line process with filtration $\F_0\vee \F^M$ satisfying (F4). If there exists a continuous regular version $P_0^{(0)}$  of the conditional avoidance probability function of $M$ given $\F_0$, then     the $(\F_0\vee\F^M)^*$-compensator of $M$ is
     \begin{equation}\label{7}
\tM^*(t)=\Lambda_0(A_t\cap \xi_1(M)),
\end{equation} where  $\Lambda_0=-\ln P_0^{(0)}$.
 Furthermore, if $Q=P|_{{\cal F}_0}$, then the law of $M$ is characterized by $Q$ and $\tM^*$.
\end{Th}
%If $\Lambda$ has a density $\lambda$, then
%\begin{eqnarray*}\tM^*(t)&=&\int _{A_t} I(M(u-)=0) \lambda(u) du\\
%&=&\Lambda(A_t\cap \xi_1(M))=-\ln P_0(A_t\cap \xi_1(M)).
%\end{eqnarray*}
%Therefore, the law of $M$ is characterized by $\tM^*$.  \end{Co}

\noindent {\bf Note:} %Suppose that $\F_0$ is trivial. 
Compare (\ref{7}) with (\ref{3}),  the formula for the compensator of   a single jump process $M$ on $\RR_+$ (with $\F_0$ trivial).  If the jump point of $M$ has continuous distribution $F$, then $P_0=1-F$   and from (\ref{3}), the compensator is $-\ln P_0(t\wedge \tau_1)=\Lambda(A_t\cap \xi_1(M))$. Thus, (\ref{7}) and (\ref{3}) are identical and in both cases, $\Lambda=-\ln P_0$ can be interpreted as a cumulative hazard.  The same will be true if $\F_0$ is not trivial.

  \section{The compensator of a general point process}\label{section7}

 We are now ready to develop a recursive formula for  the general point process compensator.  Let $N$ be a general strictly simple point process on $\RR_+^2$ with filtration $\F=\F_0\vee\F^N$ satisfying (F4) and let $N =\sum_{k=1}^\infty M_k$ be the decomposition into single line point processes of Definition \ref{4.2}.  %Let $P_0^{(k)}$ denote the conditional avoidance probability of the single line process $M_k$, given $\F(\xi_{k-1})$. 
 We will proceed as follows, letting $k\geq 1$:
    \begin{enumerate}
\item We will show that if the filtration $\F(t)=\F_0\vee \F^N(t)$ satisfies (F4) under $P$, then so does $\G(t):=\F(\xi_{k-1}(N))\vee  \F^{M_k}(t)$.   This is the key point in the development of the general point process compensator.
   \item Since $\xi_k(N)=\xi_{k-1}^+(N)\cap\xi_1(M_k)$ and $\xi_{k-1}^+(N)$ is $\F(\xi_{k-1}(N))$-measurable, the conditional law of $\xi_k(N)$  given $\F(\xi_{k-1}(N))$ is determined by the conditional law of $\xi_1(M_k)$.  By Theorem \ref{6.3} and the preceding point, this in turn is characterized by the conditional avoidance probability function 
   \begin{equation}\label{8a} P_0^{(k)}(t):=P( M_k(t)=0 \mid \F(\xi_{k-1}(N))).
   \end{equation}  Therefore, the law of $N$ is determined by $Q=P|_{{\cal F}_0}$ and the conditional avoidance probability functions $P_0^{(k)}, k\geq 1$. 
\item    %For each $t$, $P_0^{(k)}(t)I(t\in\xi_{k-1}^c)$ is $\F(A_t\cap \xi_{k-1}(N))$-measurable.
Define  $ \Lambda_k(A_t, \omega):=-\ln P_0^{(k)}(t,\omega)$. % when $t\in (\xi_{k-1}^+(N)\setminus \xi_{k-1}(N) )$. 
Letting $ \F(\xi_{k-1}(N))$ play the role of $\F_0$ for $M_k$ and defining $\G(t)$ as in point 1 above, it will be shown that   %the compensator of $M_k$ is determined by $ \Lambda_k(\cdot, \omega):=-\ln P_0^{(k)}(\cdot,\omega)$. Note that $ \Lambda_k(\cdot, \omega)$ is   a measure with support   $\xi_{k-1}^+(N)(\omega)\setminus \xi_{k-1}(N) (\omega)$ for each $\omega$. Remembering that $\xi_{k-1}^+(N)$ is $\F(\xi_{k-1}(N))$-measurable and referring to (\ref{7}),    
$$\tM^*_k(t)=\Lambda_k(A_t\cap \xi_k(N))I(t\in \xi_{k-1}^c(N)).$$
is both the $\G^*$- and the $\F^*$- compensator of $M_k$. 
%SHOW THIS  and furthermore, (SHOW THIS), for each $t$,$P_0^{(k)}(t)$ is $\F(A_t\mid \xi)$-measurable.
    \item Since $\tN^*=\sum_k\tM_k^*$,  when (F4) holds the law of $N$ is therefore characterized by $Q$ and $\tN^*$. 

% \begin{itemize}
 % \item 
 % \item<1-> $P_0^{(k)}(\omega)$ is a regular version of the conditional avoidance probability if $ \Lambda_k(\omega)=-\ln P_0^{(k)(\omega)$ is   a measure on $\xi_{k-1}^+(\omega)\setminus \xi_{k-1} (\omega)$ for a.a. $\omega$.
  \end{enumerate}
%  \begin{Df}      $P_0^{(k)}(\cdot,\omega)$ is a regular version of the conditional avoidance probability if  $ \Lambda_k(\cdot, \omega):=-\ln P_0^{(k)}(\cdot,\omega)$ is   a measure with support   $\xi_{k-1}^+(N)(\omega)\setminus \xi_{k-1}(N) (\omega)$ for each $\omega$.

 % \end{Df}

% REGULARITY ETC. START HERE
\noindent Putting the preceding points together, we arrive at  our main result:
 \begin{Th}\label{7.1} Let $N$ be a strictly simple point process on $\RR_+^2$ with filtration $(\F(t)=\F_0\vee \F^N(t))$ satisfying (F4). Assume that there exists a continuous regular version of $P_0^{(k)}~\forall k\geq 1$, where $P_0^{(k)}$ is as defined in (\ref{8a}). Then the *-compensator of $N$ has the regenerative form:
 \begin{equation}\label{8}
\tN ^*(t)=\sum_{k=1}^\infty \Lambda_k(A_t\cap \xi_k(N))I(t\in \xi_{k-1}^c(N))
 \end{equation}
 where $\Lambda_k(t) = -\ln P_0^{(k)}(t)$. % on $\{ t\in (\xi_{k-1}^+(N)\setminus \xi_{k-1}(N) )\}$.
  If $Q=P|_{{\cal F}_0}$, then the law of $N$ is characterized by $Q$ and $\tN^*$.
 \end{Th}

 \begin{Com} {\rm Theorem \ref{7.1} is the two-dimensional analogue of  the corresponding result for point processes on $\RR_+$, and in fact  the formulas in one and two dimensions are identical:  recalling (\ref{4}) (the compensator on $\RR_+$),
\begin{eqnarray*} 
\tN (t)&=&\sum_{n=1}^\infty \Lambda_n(t\wedge \tau_n)I(\tau_{n-1}<t)\\
&=&\sum_{n=1}^\infty \Lambda_n(A_t\cap \xi_n(N))I(t\in \xi_{n-1}^c(N)),
\end{eqnarray*}
%The formula for point processes on  $\RR_+$ is identical, since on $\RR_+$, $\xi_n=(0, \tau_n]$,  $ \Lambda_n(A_t\cap \xi_n)= \Lambda_n( t\wedge \tau_n)$  and $I(t\in \xi_{n-1}^c)=I(\tau_{n-1}<t)$.
 which is the same as (\ref{8}).}
  \end{Com}
  
\noindent{\bf Proof of Theorem \ref{7.1}}:\\
We must fill in the details of points 1-4, listed above.
  \begin{enumerate}
  \item \begin{itemize} \item We will begin by showing that    that for any $\F$- stopping set $\xi$ and incomparable points $s,t\in \RR_+^2$, 
  $$ \F(s) {\perp} \F(t)\; |\;({{{\cal F}(\xi)\vee{\cal F}(s\wedge t)} }),$$
  or equivalently that for any  $F\in   \F(t)$,
  \begin{equation}\label{10} P(F\mid\F(\xi)\vee\F(s))=P(F\mid\F(\xi)\vee\F(s\wedge t)).\end{equation}
This then shows that 
$$ (\F(s)\vee{\cal F}(\xi)) {\perp}( \F(t)\vee{\cal F}(\xi))\; |\;({{{\cal F}(\xi)\vee{\cal F}(s\wedge t)} } )$$
and so if $\G(s):=\F(\xi_{k-1})\vee \F^{M_k}(s)$, then $\G(s) {\perp}\G(t)\; |\;({{{\cal F}(\xi_{k-1})\vee{\cal F}(s\wedge t)} })$.  
\item We will then show that for $G\in \G(t)$,
 \begin{eqnarray} P(G\mid\F(\xi_{k-1})\vee\F(s\wedge t))&=&P(G\mid\F(\xi_{k-1})\vee\F^{M_k}(s\wedge t))\nonumber\\&=&P(G\mid \G(s\wedge t)).\label{10a}\end{eqnarray}
 Since $\G(s)\subseteq \F(\xi_{k-1})\vee\F(s)$, (\ref{10}) and (\ref{10a}) prove that $\G(s)\perp\G(t)\;|\;{\G(s\wedge t)}$.
 \end{itemize} Therefore, the proof of point 1 will be complete provided that   (\ref{10}) and (\ref{10a}) are verified.\\

 To prove (\ref{10}), we recall (\ref{add}): if  (F4) holds and if $F\in\F(t)$, then  for any lower layer $D$,  $$P[F\mid\F(D)]=P[F\mid\F(D)\cap \F(t)].
$$ Next, as is shown in \cite{IMbook}, any stopping set $\xi$ can be approximated from above by a decreasing sequence $(g_m(\xi))$ of {\em discrete} stopping sets (i.e. $g_m(\xi)$ is a stopping set taking on at most countably many values in the set of lower layers and $\xi=\cap_m g_m(\xi)$). Since $\F(\xi)=\cap_m\F(g_m(\xi))$ (\cite{IMbook}, Proposition 1.5.12), it is enough to verify (\ref{10}) for $\xi$ a discrete stopping set. Let ${\cal D}$ be a countable class of lower layers such that $\sum_{D\in {\cal D}}P(\xi = D)=1$. As noted before, for $\xi$ discrete, 
  $$\F(\xi) =\{G\in \F: G\cap \{\xi=D\}\in \F(D)~ \forall D\in {\cal D}\}$$
  and it is straightforward that $\F(\xi)=\F(D)$ on $\{\xi=D\}$. For $F\in  \F(t)$, we consider $F\cap\{t\in\xi\}$ and $F\cap \{t\in\xi^c\}$ separately. First,
   \begin{eqnarray}\lefteqn{
   P(F\cap\{t\in\xi\}\mid\F(\xi)\vee\F(s))}\nonumber\\
   &=&\sum_{D\in {\cal D}} P(F\cap\{t\in\xi\}\mid\F(\xi)\vee\F(s))I(\xi = D)\nonumber\\
   &=&\sum_{D\in {\cal D}} P(F\cap\{t\in D\}\mid\F(D)\vee\F(s))I(\xi = D)\nonumber\\
    &=&\sum_{D\in {\cal D}} I(F\cap\{t\in D\})I(\xi = D)\nonumber\\
    &=&\sum_{D\in {\cal D}} P(F\cap\{t\in D\}\mid\F(D)\vee\F(s\wedge t))I(\xi = D)\nonumber\\
    &=&\sum_{D\in {\cal D}} P(F\cap\{t\in \xi\}\mid\F(\xi)\vee\F(s\wedge t))I(\xi = D)\nonumber\\
        &=&  P(F\cap \{t\in \xi\}\mid\F(\xi)\vee\F(s\wedge t)).\label{first} 
   \end{eqnarray}
  Next,   \begin{eqnarray}\lefteqn
  { P(F\cap\{t\in\xi^c\}\mid\F(\xi)\vee\F(s))}\nonumber\\&=&\sum_{D\in {\cal D}} P(F\cap\{t\in\xi^c\}\mid\F(\xi)\vee\F(s))I(\xi = D)\nonumber\\
   &=&\sum_{D\in {\cal D}} P(F\cap\{t\in D^c\}\mid\F(D)\vee\F(s))I(\xi = D)\nonumber\\
   &=&\sum_{D\in {\cal D}} P(F\cap\{t\in D^c\}\mid\F(D\cup A_s))I(\xi = D)\nonumber\\
  &=&\sum_{D\in {\cal D}} P(F\cap\{t\in D^c\}\mid\F(D\cup A_s)\cap \F(t))I(\xi = D)\label{11}\\
 &=&\sum_{D\in {\cal D}} P(F\cap\{t\in D^c\}\mid\F((D\cup A_s)\cap A_t)I(\xi = D) \nonumber\\
     &=&\sum_{D\in {\cal D}} P(F\cap\{t\in D^c\}\mid\F(D\cup (A_s\cap A_t))I(\xi = D) \label{12}\\
      &=&\sum_{D\in {\cal D}} P(F\cap\{t\in D^c\}\mid\F(D)\vee \F(A_s\cap A_t))I(\xi = D) \nonumber\\
       &=&\sum_{D\in {\cal D}} P(F\cap\{t\in\xi^c\}\mid\F(\xi)\vee \F(A_s\cap A_t))I(\xi = D) \nonumber\\
        &=&  P(F\cap\{t\in\xi^c\}\mid\F(\xi)\vee\F(s\wedge t)). \label{second} 
   \end{eqnarray}
  Equations (\ref{11}) and (\ref{12}) follow from (\ref{add}). Putting (\ref{first}) and (\ref{second}) together yields (\ref{10}).
  
%  Now we prove (\ref{10a}).  Since $s$ and $t$ are incomparable, without loss of generality we will assume that $t_1<s_1$ and $t_2>s_2$ and so $s\wedge t=(t_1,s_2)$. We have $\F(\xi_{k-1})\vee\F(s\wedge t)=\F(\xi_{k-1})\vee\F^N(s\wedge t)$, and so on  the event $\{M_k(s\wedge t)=0\}$, $\F(\xi_{k-1})\vee\F(s\wedge t)=\F(\xi_{k-1})\vee\F^{M_k}(s\wedge t)$  because there are no jumps in $A_{s\wedge t}\setminus \xi_{k-1}$. Therefore, (\ref{10a}) is satisfied on $\{M_k(s\wedge t)=0\}$. 
  
%  To handle the event $\{M_k(s\wedge t)>0\}$, let $\tau=\inf(v\in \RR_+:M_k(v,s_2)>0)\wedge t_1$.  We have that $\tau$ is a stopping time with respect to the one-dimensional filtration $\F(\xi_{k-1})\vee\F^{M_k}(\cdot,s_2)$ and $\tau\leq t_1$ on $\{M_k(s\wedge t)>0\}$. Note that $\F(\xi_{k-1})\vee\F(\tau,s_2)=\F(\xi_{k-1})\vee\F^{M_k}(\tau,s_2)$ since there are no jumps on $A_{(\tau,s_2)}\setminus \xi_{k-1}$ other than (possibly) a single jump from $M_k$ on the line segment $\{(\tau,u), 0\leq u\leq s_2\}$. Approximate $\tau$ from above with discrete stopping times $\tau_m\leq t_1$, $\tau_m\downarrow \tau$. By right continuity of the filtrations,

  Now we prove (\ref{10a}).  Since $s$ and $t$ are incomparable, without loss of generality we will assume that $t_1<s_1$ and $t_2>s_2$ and so $s\wedge t=(t_1,s_2)$. We have $\F(\xi_{k-1})\vee\F(s\wedge t)=\F(\xi_{k-1})\vee\F^N(s\wedge t)$. Let $\tau:=\inf(v\in \RR_+:M_k(v,s_2)>0)\wedge t_1$;  $\tau$ is a stopping time with respect to the one-dimensional filtration $\F(\xi_{k-1})\vee\F^{M_k}(\cdot,s_2)$.  Note that $\F(\xi_{k-1})\vee\F(\tau,s_2)=\F(\xi_{k-1})\vee\F^N(\tau,s_2)=\F(\xi_{k-1})\vee\F^{M_k}(\tau,s_2)$ since $N$ has  no jumps on $A_{(\tau,s_2)}\setminus \xi_{k-1}$ other than (possibly) a single jump from $M_k$ on the line segment $\{(\tau,u), 0\leq u\leq s_2\}$. Approximate $\tau$ from above with discrete stopping times $\tau_m\leq t_1$, $\tau_m\downarrow \tau$. By right continuity of the filtrations,
 \begin{eqnarray}
  \F(\xi_{k-1})\vee\F(\tau_m,s_2) &=&\F(\xi_{k-1})\vee\F^N(\tau_m,s_2)\nonumber\\
  &\downarrow &\F(\xi_{k-1})\vee\F^N(\tau,s_2)\nonumber\\
 & =&\F(\xi_{k-1})\vee\F^{M_k}(\tau,s_2). \label{rt}
  \end{eqnarray}
 
 Without loss of generality, let $G=\{M_k(t)=j\}$ in (\ref{10a}). Observe that  on $\{M_k(s\wedge t)>0\}$, 
$M_k(t)=M_k(\tau_m,t_2)+M_k(((\tau_m,0),(t_1,s_2)])$
for every $m$, since $\{M_k(\tau_m)>0\}$ and the jumps of $M_k$ are incomparable. On $\{M_k(s\wedge t)=0\}$, $\tau_m=t_1\forall  m$ and $M_k(t)=M_k(\tau_m,t_2)$. For ease of notation in what follows, let $X(\tau_m):=M_k(\tau_m,t_2)$ and $Y(\tau_m):=M_k((\tau_m,0),(t_1,s_2)])$. Recall that $s\wedge t=(t_1,s_2)$ and let $R_m$ denote the (countable) set of possible values of $\tau_m$.
\begin{eqnarray}\lefteqn{
P(M_k=j\mid\F(\xi_{k-1})\vee\F(s\wedge t))}\nonumber\\
 &=&\sum_{r\in R_m}P(M_k=j\mid\F(\xi_{k-1})\vee\F(s\wedge t))I(\tau_m=r)\nonumber\\
&=&\sum_h\sum_{r\in R_m}P(X(\tau_m)=h,Y(\tau_m)=j-h\mid\F(\xi_{k-1})\vee\F(s\wedge t))I(\tau_m=r)\nonumber\\
&=&\sum_h\sum_{r\in R_m}P(X(r)=h\mid\F(\xi_{k-1})\vee\F(s\wedge t))I(Y(r)=j-h)I(\tau_m=r)\nonumber\\
%&=&\sum_h\sum_rP(X(r)=h\mid\F(\xi_{k-1})\vee\F(r,s_2))I((Y(r)=j-h)I(M_k(s\wedge t)>0)I(\tau_m=r)\lnonumber\\
&=&\sum_h\sum_{r\in R_m}P(X(r)=h\mid\F(\xi_{k-1})\vee\F(r,s_2))\nonumber\\
&&~~~~~~~~~~~~~~~~~~~~~~~\times I(Y(r)=j-h)I(\tau_m=r)\label{17b}\\
&=&\sum_h\sum_{r\in R_m}P(X(r)=h,Y(r)=j-h\mid\F(\xi_{k-1})\vee\F(r,s_2)\vee \F^{M_k}(s\wedge t))\nonumber\\
&&~~~~~~~~~~~~~~~~~~~~~~~\times I(\tau_m=r)\label{17b1}\\
&=&P(M_k=j\mid\F(\xi_{k-1})\vee\F(\tau_m,s_2)\vee \F^{M_k}(s\wedge t))\nonumber\\
&\stackrel{{m\rightarrow\infty}}{\longrightarrow}&P(M_k=j\mid\F(\xi_{k-1})\vee\F^{M_k}(\tau,s_2)\vee \F^{M_k}(s\wedge t))\label{17c}\\
&=&P(M_k=j\mid\F(\xi_{k-1})\vee \F^{M_k}(s\wedge t)).\nonumber
\end{eqnarray}
%(\ref{17a}) and (\ref{17d}) follow since $r< t_1$ on $\{M_k(s\wedge t)>0\}$, 
(\ref{17b}) and (\ref{17b1}) follow from (\ref{add}) and the fact that $X(r)$ is $\F(r,t_2)$-measurable, and (\ref{17c}) follows from (\ref{rt}). %Taking the events $\{M_k(s\wedge t)=0\}$ and $\{M_k(s\wedge t)>0\}$ 
This proves (\ref{10a}) and completes the proof of point 1.
   \item This follows immediately from point 1 and Theorem \ref{6.3}.
   \item %The argument follows closely that used in one dimension (see \cite{DVJ}, pp. 520-521).
    Begin by recalling that $M_k$ has its support on $\xi_{k-1}^+(N)\setminus \xi_{k-1}(N)$
     and  so $P_0^{(k)}(t)=P(M_k(t)=0|\F(\xi_{k-1}(N)))=P(M_k(A_t \cap \xi_{k-1}^+(N))=0|\F(\xi_{k-1}(N)))$. %on $\{ t\in(\xi_{k-1}^+(N))^c\}$. Likewise, on $\{t\in \xi_{k-1}(N)\}$,  $P_0^{(k)}(t)=1$. 
    Therefore, we will identify $M_k$ with the first line of a Cox process whose driving measure $\Lambda_k (t)=-\ln  P_0^{(k)}(t)$ has support $\xi_{k-1}^+(N)\setminus \xi_{k-1}(N)$.
   % $\Lambda_k$ as follows:  on $\{t\in\xi_{k-1}^+(N)\} $, 
    %We now have a well-defined measure on $\xi_{k-1}^+(N)$ and on  $\{t\in (\xi_{k-1}^+(N) )^c\}$, $ \Lambda_k(t) =\Lambda (A_t\cap \xi_{k-1}^+(N))$. % $$\xi_{k-1}^+(N)$$ \Lambda_k(t):=-\ln( P_0^{(k)}(A_t\cap(\xi_{k-1}^+(N)\setminus \xi_{k-1}(N) )).$$  
  Now, identifying $\G_0=\F(\xi_{k-1}(N))$ and $\G(t)=\F(\xi_{k-1}(N))\vee \F^{M_k}(t)$,    as in Theorem \ref{6.4}  we have  that the $\G^*$-compensator of $M_k$ is:
 \begin{eqnarray}
\tM^*_k(t)&=&\Lambda_k(A_t\cap\xi_1(M_k))\nonumber \\
&=&\Lambda_k(A_t\cap\xi_1(M_k))I(t\in \xi_{k-1}^c(N)  )\nonumber \\
&=&\Lambda_k(A_t\cap  \xi_{k-1}^+(N) \cap \xi_1(M_k))I(t\in \xi_{k-1}^c(N)   )\nonumber\\
%&=&\Lambda_k(A_t\cap \xi_{k} (N)\setminus \xi_{k-1}(N) ))I(t\in \xi_{k-1}^c(N)  ) \label{15}\\
%&=&-\ln( P_0^{(k)}(A_t\cap \xi_{k} (N)) I(t\in  \xi_{k-1}^c(N) )\nonumber\\
&=&\Lambda_k(A_t\cap \xi_{k} (N)) I(t\in \xi_{k-1}^c(N)  ).\nonumber
 \end{eqnarray}  
The last two equalities follow  since  $\{ t\in\xi_{k-1}^+(N)\}$ is   $\F(\xi_{k-1}(N))$-measurable and $\xi_{k}(N)=\xi_{k-1}^+(N)\cap \xi_1(M_k)$. 

We must now show that $\tM^*_k$ is also the $\F^*$-compensator of $M_k$. First we show that $\tM^*_k$ is $\F$-adapted. On $\{t\in \xi_{k-1}\}\in \F(t)$, $P_0^{(k)}(t)=0$. On $\{t\in \xi_{k-1}^c\}$,    by (\ref{add}) and taking discrete approximations of $\xi_{k-1}$, arguing as in the proof of (\ref{second}) we have 
$$P_0^{(k)}(t)I(t\in \xi_{k-1}^c)=P(M_k(t)=0\mid \F(\xi_{k-1})\cap \F(t))I(t\in \xi_{k-1}^c).$$
 Therefore, $-\ln P_0^{(k)}$ is $\F$-adapted. Since $\tM^*_k$ is $\F$-adapted and continuous, by Theorem \ref{strongmg} it remains only to prove that  
 $$E\left[(M_k-\tM_k^*)(s,t]\mid\F^*(s)\right]=0.$$ 
 First,  if  $t\in \xi_{k-1}$ or if  $ s\in (\xi_{k-1}^+)^c $, then $(M_k-\tM_k^*)(s,t]=0$ and so trivially
    %Next, note that $\{s\in (\xi_{k-1}^+)^c\}$ is  $\F^*(s)$-measurable. Therefore,
     \begin{equation}\label{16}
    E\left[(M_k-\tM_k^*)(s,t]I(t\in  \xi_{k-1})\mid\F^*(s)\right]=0,
    \end{equation}
   and
    \begin{equation}\label{17}
    E\left[(M_k-\tM_k^*)(s,t]I(s\in (\xi_{k-1}^+)^c)\mid\F^*(s)\right]=0.
    \end{equation}   
   For $s<t$, $t\in \xi_{k-1}^c$ and $ s\in  \xi_{k-1}^+ $, it is enough to show that
  \begin{eqnarray} \lefteqn{E\left[ M_k (s,t]I( s\in  \xi_{k-1}^+,t\in \xi_{k-1}^c)\mid \G^*(s)\right]}\nonumber\\ 
&=&E\left[M_k (s,t]I( s\in  \xi_{k-1}^+,t\in \xi_{k-1}^c) \mid\F(\xi_{k-1})\vee (\F^{M_k})^*(s)\right]\nonumber\\
&=&E\left[M_k (s,t] I( s\in  \xi_{k-1}^+,t\in \xi_{k-1}^c)\mid\F(\xi_{k-1})\vee (\F^{N})^*(s)\right].\label{crucial} 
\end{eqnarray} 
If (\ref{crucial}) is true, then since $\F^*(s)=\F_0\vee (\F^{N})^*(s)\subseteq \F(\xi_{k-1})\vee (\F^{N})^*(s)$,
\begin{eqnarray} 0&=& E\left[(M_k-\tM_k^*)(s,t]I( s\in  \xi_{k-1}^+,t\in \xi_{k-1}^c)\mid \G^*(s)\right] \nonumber\\
%&=&E\left[(M_k-\tM_k^*)(s,t]\mid\F(\xi_{k-1})\vee (\F^{M_k})^*(s)\right]\nonumber\\
&=&E\left[(M_k-\tM_k^*)(s,t]I( s\in  \xi_{k-1}^+,t\in \xi_{k-1}^c)\mid\F(\xi_{k-1})\vee (\F^{N})^*(s)\right]\nonumber\\
&=&E\left[(M_k-\tM_k^*)(s,t]I( s\in  \xi_{k-1}^+,t\in \xi_{k-1}^c)\mid\F^*(s)\right].\label{crucial2}
\end{eqnarray}
 To prove (\ref{crucial}), let $\tau_1=\inf\{v:M_k(v,s_2)>0\}\wedge t_1$. Similar to the  argument   used to prove   (\ref{10a}), 
  we have
  $$\F(\xi_{k-1})\vee \F^N(\tau_1,s_2)=\F(\xi_{k-1})\vee \F^{M_k}(\tau_1,s_2)$$
  and  using (F4) (cf. (\ref{add})) and discrete approximations for $\tau_1$, it follows that
  \begin{eqnarray*}\lefteqn{E\left[M_k (s,t] I( s\in  \xi_{k-1}^+,t\in \xi_{k-1}^c)\mid\F(\xi_{k-1})\vee (\F^{N})^*(s)\right]}\\
  &=&E\left[M_k (s,t] I( s\in  \xi_{k-1}^+,t\in \xi_{k-1}^c)\mid\F(\xi_{k-1})\vee (\F^{N})^1(s)\vee\F^{M_k}(t_1,s_2)\right].\end{eqnarray*}
Next, letting $\tau_2=\inf\{u:M_k(s_1,u)>0\}\wedge t_2$, we argue as above and apply (F4) (cf. (\ref{add})) twice to obtain
\begin{eqnarray*}\lefteqn
{ E\left[M_k (s,t] I( s\in  \xi_{k-1}^+,t\in \xi_{k-1}^c)\mid\F(\xi_{k-1})\vee (\F^{N})^1(s)\vee\F^{M_k}(t_1,s_2)\right]}\nonumber\\
&=&E\left[M_k (s,t] I( s\in  \xi_{k-1}^+,t\in \xi_{k-1}^c)\mid\F(\xi_{k-1})\vee \F^{M_k}(s_1,t_2)\vee\F^{M_k}(t_1,s_2)\right]\\
&=&E\left[M_k (s,t] I( s\in  \xi_{k-1}^+,t\in \xi_{k-1}^c)\mid\G^*(s)\right].\end{eqnarray*} 
This completes the proof of (\ref{crucial}) and (\ref{crucial2}). Combining (\ref{16}), (\ref{17}) and (\ref{crucial2}), it follows that $\tM_k^*$ is the $\F^*$-compensator of $M_k$.

\item This is immediate because of the decomposition $N=\sum_{k=1}^\infty M_k$.
  \end{enumerate} This completes the proof of Theorem \ref{7.1}.\hfill $\Box$

  \section{Conclusion}\label{section8}
  In this paper we have proven a two-dimensional analogue of Jacod's characterization of the law of a point process via a regenerative formula for its compensator. For clarity we have restricted our attention to continuous avoidance probabilities. There remain  many open questions that merit  further investigation, for example:
  \begin{itemize}
%  \item add an initial $\sigma$-field
  \item Extend the regenerative formula to discontinuous avoidance probability functions. In this case, the logarithmic relation between the avoidance probability and the cumulative hazard will be replaced by a product limit formula. 
  \item Extend the regenerative formula to marked point processes.
    \item Find a complete characterization of the class of predictable increasing functions that are *-compensators for  planar point processes satisfying (F4), in analogy to Theorem 3.6 of \cite{Jacod}.
  \item Generalize the results of this paper to point processes on $\RR_+^d,d>2$. The main challenge will be to find an appropriate $d$-dimensional analogue of (F4).
   \end{itemize}

\end{document}